\documentclass[a4paper,12pt]{article}
\usepackage[left=2.50cm, right=2.50cm, top=1cm, bottom=2.50cm]{geometry} 
\usepackage{amsfonts,amsmath,latexsym,amssymb,amsthm,mathrsfs,upref} 
\usepackage{enumitem}        
\usepackage{caption}

\newtheorem{lem}{Lemma}
\newtheorem{thm}{Theorem}

\newtheorem{definition}{Definition}

\newcommand\keywords[1]{\quad\quad\, \textbf{Keywords}: #1}

\title{The classification of $n$-dimensional nilpotent non-Tortkara anticommutative algebras with $(n-4)$-dimensional annihilator}%
\author{
	\small {Shun Xu}\\
	{\small School of Mathematical Sciences, Tongji University, Shanghai, 200092, China}\\
	{\small  e-mail: shunxu1997@163.com}
}
\date{}

\begin{document}
	\maketitle
	\begin{abstract}
		In this paper, we give a complete classification of $n$-dimensional nilpotent non-Tortkara anticommutative algebras with $(n-4)$-dimensional annihilator over $\mathbb{C}$.
	\end{abstract}
	
\section{Introduction}

One of the classical problems in the theory of  algebras is to classify (up to isomorphism) the nilpotent algebras of dimension $n$ from a certain variety defined by some family of polynomial equalities. There
are many results related to algebraic classification of small dimensional nilpotent algebras in varieties of Jordan,
Lie, Leibniz, Zinbiel and other algebras. Algebraic classification of nilpotent Lie algebras of dimension 7 (over algebraically closed field and $\mathbb{R}$)  \cite{gong1998classification}; algebraic classification of five-dimensional nilpotent Jordan algebras \cite{hegazi2016classification2}; algebraic classification of four-dimensional nilpotent Leibniz algebras \cite{demir2017classification};  algebraic  classification of complex 5-dimensional Zinbiel algebras \cite{alvarez2022algebraic}.

There are also some results related to algebraic classification of $n$-dimensional nilpotent algebras with $(n-s)$-dimensional annihilator, where $s$ is small positive integer. The classification of $n$-dimensional anticommutative algebras with $(n-3)$-dimensional annihilator \cite{calderon2019classification}; The classification of $n$-dimensional non-associative Jordan algebras with $(n-3)$-dimensional annihilator \cite{hegazi2018classification};  The classification of $n$-dimensional non-Lie Malcev algebras with $(n-4)$-dimensional annihilator \cite{hegazi2016classification}.

 In this paper, we give a complete classification of $n$-dimensional nilpotent non-Tortkara anticommutative algebras with $(n-4)$-dimensional annihilator over $\mathbb{C}$. An  algebra $(A,\cdot)$ over  $\mathbb{C}$ is anticommutative if it satisfies: $xy+yx=0$ for all $x, y \in A$. Let $A$ be an anticommutative algebra. The ideal $\mathrm{A n n}(A)=\{x \in A:xA=0\}$ is called the annihilator of $A$. An anticommutative
algebra is Tortkara algebra, if it satisfies:
$$(ab)(cb) = J(a, b, c)b, \text{where}\, J(a, b, c)=(ab)c + (bc)a + (ca)b.$$

Our main result in this paper is the following theorem.

\begin{thm}[Main result]\label{thm1}
	Let $A$ be a $n$-dimensional nilpotent non-Tortkara anticommutative algebra  with $\dim \mathrm{Ann}(A)=n-4$ 
	over $\mathbb{C}$. Then $A$ is isomorphic to one of the followings:
	\begin{itemize}
		\item 
		$n=4$
		\begin{itemize}
			\item $\mathrm{empty\,\,set};$
		\end{itemize}
		
		\item
		$n=5$
		\begin{itemize}
			\item 
			$A_{5,1}:e_{1} e_{2}=e_{3}, e_{1} e_{3}=e_{4}, e_{3} e_{4}=e_{5};$
		\end{itemize}
		\item
		$n=6$
		\begin{itemize}
			\item 
			$A_{6,1}=A_{5,1}\oplus \mathbb{C}e_6;$
			\item
			$A_{6,2}:e_{1} e_{2}=e_{3}, e_{1} e_{3}=e_{4}, e_{3} e_{4}=e_{5}, e_{2} e_{4}=e_{6};$
			\item
			$A_{6,3}:e_{1} e_{2}=e_{3}, e_{1} e_{3}=e_{4}, e_{3} e_{4}=e_{5}, e_{2} e_{3}=e_{6};$
			\item
			$A_{6,4}:e_{1} e_{2}=e_{3}, e_{1} e_{3}=e_{4}, e_{3} e_{4}=e_{5}, e_{1} e_{4}=e_{6};$
			\item
			$A_{6,5}:e_{1} e_{2}=e_{3}, e_{1} e_{3}=e_{4}, e_{3} e_{4}=e_{5}, e_{2} e_{3}=e_{6}, e_{1} e_{4}=e_{6};$
		\end{itemize}
		
		\item
		$n=7$
		\begin{itemize}
			\item 
			$A_{7,i}=A_{6,i}\oplus\mathbb{C}e_7,i=1,\cdots,5;$
			\item
			$A_{7,6}:e_{1} e_{2}=e_{3}, e_{1} e_{3}=e_{4}, e_{3} e_{4}=e_{5}, e_{2} e_{4}=e_{6}, e_{2} e_{3}=e_{7};$
			\item
			$A_{7,7}:e_{1} e_{2}=e_{3}, e_{1} e_{3}=e_{4}, e_{3} e_{4}=e_{5}, e_{2} e_{4}=e_{6}, e_{1} e_{4}=e_{7};$
			\item
			$A_{7,8}:e_{1} e_{2}=e_{3}, e_{1} e_{3}=e_{4}, e_{3} e_{4}=e_{5}, e_{1} e_{4}=e_{6}, e_{2} e_{3}=e_{7};$
			\item
			$A_{7,9}:e_{1} e_{2}=e_{3}, e_{1} e_{3}=e_{4}, e_{3} e_{4}=e_{5}, e_{2} e_{4}=e_{6}, e_{1} e_{4}=e_{7}, e_{2} e_{3}=e_{7};$
		\end{itemize}
		
		\item
		$n=8$
		\begin{itemize}
			\item 
			$A_{8,i}=A_{7,i}\oplus\mathbb{C}e_8,i=1,\cdots,9;$
			\item
			$A_{8,10}:e_{1} e_{2}=e_{3}, e_{1} e_{3}=e_{4}, e_{3} e_{4}=e_{5}, e_{2} e_{4}=e_{6}, e_{2} e_{3}=e_{7},e_{1} e_{4}=e_{8};$
		\end{itemize}
		
		\item
		$n\geqslant9$
		\begin{itemize}
			\item 
			$
			A_{n, i}=A_{8, i} \oplus \mathbb{C} e_{9} \oplus \cdots \oplus \mathbb{C} e_{n}, i=1, \ldots, 10.
			$
		\end{itemize}
	\end{itemize}
\end{thm}
The paper is organized as follows. In Section \ref{sec2}, we describe a method(Skjelbred-Sund method) for classifying all anticommutative algebras of dimension $n$ with $s$-dimensional annihilator given those algebras of dimension $n-s$, which also
appeared in \cite{hegazi2016classification}. In Section \ref{sec3}, the proof of Theorem \ref{thm1} is given.

Throughout the paper we use the following notation. All the vector spaces and algebras will be assumed to be finite dimensional over $\mathbb{C}$. The multiplication of an algebra is specified by giving only the nonzero products among the basis elements.

\section{The analogue of Skjelbred-Sund method for anticommutative algebras}\label{sec2} 

 Hegazi given the analogue of the Skjelbred–Sund method for Malcev algebras in \cite{hegazi2016classification}. Inspired by \cite{hegazi2016classification}, we give the analogue of the Skjelbred–Sund method for anticommutative algebras in this section. Proofs for all the conclusions mentioned in this section can be found in \cite{hegazi2016classification}.

Let ${A}$ be an anticommutative algebra, ${V}$ be a vector space and ${Z}^{2}({A}, {V})$ is defined to be the set of all skew-symmetric bilinear maps $\theta: {A} \times {A} \longrightarrow {V}$. For $f \in \mathrm{Hom}({A}, {V})$, we define $\delta f \in {Z}^{2}({A}, {V})$ by the equality $\delta f(x, y)=f(x y)$ and set ${B}^{2}({A}, {V})=\{\delta f \mid f \in \operatorname{Hom}({A}, {V})\}$. One can easily check that ${B}^{2}({A}, {V})$ is a linear subspace of ${Z}^{2}({A}, {V})$. Let us define ${H}^{2}({A}, {V})$ as the quotient space ${Z}^{2}({A}, {V}) / {B}^{2}({A}, {V})$. The equivalence class of $\theta \in {Z}^{2}({A}, {V})$ in ${H}^{2}({A}, {V})$ is denoted by $[\theta]$. As usual, we call the elements of ${Z}^{2}({A}, {V})$ cocycles, those of ${B}^{2}({A}, {V})$ coboundaries, and ${H}^{2}({A}, {V})$ is the corresponding second cohomology space.

Suppose now that $\operatorname{dim} {A}=m<n$ and $\operatorname{dim} {V}=n-m$. Let $\theta \in Z^{2}(A, V)$, we can define on the space ${A}_{\theta}:={A} \oplus {V}$ the anticommutative bilinear product  by the equality $(x+x^{\prime})(y+y^{\prime})=x y+\theta(x, y)$ for $x, y \in {A}, x^{\prime}, y^{\prime} \in {V}$. $A_{\theta}$ is an anticommutative algebra and  ${A}_{\theta}$ is nilpotent if and only if ${A}$ is nilpotent. The set $\theta^{\perp}=\{x \in A: \theta(x, A)=0\}$ is called the radical of $\theta$. Then  $\mathrm{Ann}\left(A_{\theta}\right)=\left(\theta^{\perp} \cap \mathrm{A n n}(A)\right) \oplus V$ by \cite[Lemma 4]{hegazi2016classification}. The algebra ${A}_{\theta}$ is called a $(n-m)$-dimensional central extension of ${A}$ by ${V}$.

\begin{lem}\cite[Lemma 5]{hegazi2016classification}
	Let $A$ be an anticommutative algebra with  $\mathrm{Ann}(A) \neq 0$. Then there exists, up to isomorphism, a unique anticommutative algebra $A^{\prime}$, and $\theta \in Z^{2}\left(A^{\prime}, \mathrm{Ann} (A)\right)$ with $\theta^{\perp} \cap \mathrm{Ann} \left(A^{\prime}\right)=0$ such that $A \cong A_{\theta}^{\prime}$ and $A / \mathrm{A n n}(A) \cong A^{\prime} .$
\end{lem}

Let $e_{1}, \ldots, e_{s}$ be a basis of $V$, and $\theta \in Z ^{2}(A, V)$. Then $\theta$ can be uniquely written as $\theta(x, y)=\sum_{i=1}^{s} \theta_{i}(x, y) e_{i}$, where $\theta_{i} \in Z ^{2}(A, \mathbb{C}) .$ Moreover, $\theta^{\perp}=\theta_{1}^{\perp} \cap \theta_{2}^{\perp} \cap \cdots \cap \theta_{s} ^{\perp}.$ Further, $\theta \in B^{2}(A, V)$ if and only if all $\theta_{i} \in B^{2}(A, \mathbb{C}) .$

Let $A$ be an anticommutative algebra with a basis $e_{1}, e_{2}, \ldots, e_{n}$. Then by $\Delta_{i j}$ we denote the skew-symmetric bilinear form $\Delta_{i j}: A \times A \longrightarrow \mathbb{C}$ with $\Delta_{i j}\left(e_{i}, e_{j}\right)=-\Delta_{i j}\left(e_{j}, e_{i}\right)=1$ and $\Delta_{i j}\left(e_{l}, e_{m}\right)=0$ if $\{i, j\} \neq\{l, m\} .$ Then the set $\left\{\Delta_{i j}: 1 \leq i<j \leq n\right\}$ is a basis for the linear space of skew-symmetric bilinear forms on $A$. Then every $\theta \in Z ^{2}(A, \mathbb{C})$ can be uniquely written as $\theta=\sum_{1 \leq i<j \leq n} c_{i j} \Delta_{i, j}$, where $c_{i j} \in \mathbb{C} .$ 

Let $\left\{e_{1}, e_{2}, \ldots, e_{m}\right\}$ be a basis of $A^2$. Then the set $\left\{\delta e_{1}^{*}, \delta e_{2}^{*}, \ldots, \delta e_{m}^{*}\right\}$ by \cite[Lemma 6]{hegazi2016classification}, where $e_{i}^{*}\left(e_{j}\right)=\delta_{i j}$ and $\delta_{i j}$ is the Kronecker delta, is a basis of $B^{2}(A, \mathbb{C})$. Let $\theta, \vartheta \in Z^{2}(A, V)$ such that $[\theta]=[\vartheta] .$ Then $\theta^{\perp} \cap \mathrm{Ann}(A)=\vartheta^{\perp} \cap \mathrm{Ann} (A)$ or, equivalently, $\mathrm{Ann}\left(A_{\theta}\right)=\operatorname{\mathrm{\mathrm{Ann}}}\left(A_{\vartheta}\right)$ by \cite[Lemma 7]{hegazi2016classification}. Furthermore, $A_{\theta} \cong A_{\vartheta}$.

Let $\mathrm{A u t}(A)$ be the automorphism group of an anticommutative algebra $A .$ Let $\phi \in \mathrm{A u t}(A)$. For $\theta \in Z^{2}(A, V)$ define $\phi \theta(x, y)=\theta(\phi(x), \phi(y)) .$ Then $\phi \theta \in Z^{2}(A, V) .$ So, $\mathrm{Aut} (A)$ acts on $Z^{2}(A, V)$. $\phi \theta \in B^{2}(A, V)$ if and only if $\theta \in B^{2}(A, V)$ by \cite[Lemma 8]{hegazi2016classification}. So, $\mathrm{Aut} (A)$ acts on $H^{2}(A, V)$.

Let $\phi=\left(a_{i j}\right) \in A u t(A)$ and $\theta \in Z^{2}(A, \mathbb{C}) .$ Let $C=\left(c_{i j}\right)$ be the matrix representing $\theta$ and $C^{\prime}=\left(c_{i j}^{\prime}\right)$ be the matrix representing $\phi \theta$. Then $C^{\prime}=\phi^{t} C \phi$.

\begin{definition}
	Let $A$ be an anticommutative algebra. If $A=I \oplus \mathbb{C} x$ is a direct sum of two ideals, then $\mathbb{C} x$ is called an annihilator component of $A .$
\end{definition}
 Let $\theta(x, y)=\sum_{i=1}^{s} \theta_{i}(x, y) e_{i} \in Z^{2}(A, V)$ and $\theta^{\perp} \cap \mathrm{Ann}(A)=0 .$ Then $A_{\theta}$ has an annihilator component if and only if $\left[\theta_{1}\right],\left[\theta_{2}\right], \ldots,\left[\theta_{s}\right]$ are linearly dependent in $H^{2}(A, \mathbb{C})$  by \cite[Lemma 13]{hegazi2016classification}. Let $\vartheta(x, y)=\sum_{i=1}^{s} \vartheta_{i}(x, y) e_{i}$ be another element of $Z^{2}(A, V) .$ Suppose that $A_{\theta}$ has no annihilator components and $\theta^{\perp} \cap \mathrm{Ann}(A)=\vartheta^\perp \cap \mathrm{Ann}(A)=0$. Then $A_{\theta} \cong A_{\vartheta}$ if and only if there exists a map $\phi \in \mathrm{A u t}(A)$ such that the set $\left\{\left[\phi \vartheta_{i}\right]: i=1, \ldots, s\right\}$ spans the same subspace of $H^{2}(A, \mathbb{C})$ as the set $\left\{\left[\theta_{i}\right]: i=1, \ldots, s\right\} $ by \cite[Lemma 14]{hegazi2016classification}.

Let $V$ be a finite-dimensional vector space over a $\mathbb{C}$. The Grassmannian $G_{k}(V)$ is the set of all $k$-dimensional linear subspaces of $V$. Let $G_{s}\left(H^{2}(A, \mathbb{C})\right)$ be the Grassmannian of subspaces of dimension $s$ in $H^{2}(A, \mathbb{C})$. There is a natural action of $\mathrm{Aut}(A)$ on $G_{s}\left(H^{2}(A, \mathbb{C})\right)$. Let $\phi \in \mathrm{Aut} (A)$. For $W=\left\langle\left[\theta_{1}\right],\left[\theta_{2}\right], \ldots,\left[\theta_{s}\right]\right\rangle \in G_{s}\left(H^{2}(A, \mathbb{C})\right)$ define $\phi W=\left\langle\left[\phi \theta_{1}\right],\left[\phi \theta_{2}\right], \ldots,\left[\phi \theta_{s}\right]\right\rangle$. Then $\phi W \in G_{s}\left(H^{2}(A, \mathbb{C})\right)$.
We denote the orbit of $W \in G_{s}\left(H^{2}(A, \mathbb{C})\right)$ under the action of $\operatorname{\mathrm{Aut}}(A)$ by
$\mathrm{Orb}(W)$.

Let $W_{1}=\left\langle\left[\theta_{1}\right],\left[\theta_{2}\right], \ldots,\left[\theta_{s}\right]\right\rangle, W_{2}=\left\langle\left[\vartheta_{1}\right],\left[\vartheta_{2}\right], \ldots,\left[\vartheta_{s}\right]\right\rangle \in G_{s}\left(H^{2}(A, \mathbb{C})\right)$. If $W_{1}=W_{2}$, then $\bigcap_{i=1}^{s} \theta_{i}^{\perp} \cap \mathrm{A n n}(A)=\bigcap_{i=1}^{s} \vartheta_{i}^{\perp} \cap \mathrm{A n n}(A)$ by \cite[Lemma 15]{hegazi2016classification}. This result allows us to define
$$
T_{s}(A)=\left\{W=\left\langle\left[\theta_{1}\right],\left[\theta_{2}\right], \ldots,\left[\theta_{s}\right]\right\rangle \in G_{s}\left(H^{2}(A, \mathbb{C})\right): \underset{i=1}{\bigcap} \theta_{i}^{\perp} \cap A n n(A)=0\right\}
$$

The set $T_{s}(A)$ is stable under the action of  $\mathrm{Aut}(A)$  by \cite[Lemma 16]{hegazi2016classification}. 

Let $V$ be an $s$-dimensional vector space spanned by $e_{1}, e_{2}, \ldots, e_{s}$. Given a anticommutative algebra $A$, let $E(A, V)$ denote the set of all anticommutative algebras without annihilator components which are $s$-dimensional annihilator extensions of $A$ by $V$ and have $s$-dimensional annihilator. Then $E(A, V)=$ $\left\{A_{\theta}: \theta(x, y)=\sum_{i=1}^{s} \theta_{i}(x, y) e_{i}\right.$ and $\left.\left\langle\left[\theta_{1}\right],\left[\theta_{2}\right], \ldots,\left[\theta_{s}\right]\right\rangle \in T_{s}(A)\right\} .$ Given $A_{\theta} \in E(A, V)$, let $\left[A_{\theta}\right]$ denote the isomorphism class of $A_{\theta} .$ Let $A_{\theta}, A_{\vartheta} \in E(A, V)$. Suppose that $\theta(x, y)=\sum_{i=1}^{s} \theta_{i}(x, y) e_{i}$ and $\vartheta(x, y)=\sum_{i=1}^{s} \vartheta_{i}(x, y) e_{i} .$ Then $\left[A_{\theta}\right]=\left[A_{\vartheta}\right]$ if and only if $$\mathrm{Orb}\left\langle\left[\theta_{1}\right],\left[\theta_{2}\right], \ldots,\left[\theta_{s}\right]\right\rangle= \mathrm{Orb}\left\langle\left[\vartheta_{1}\right],\left[\vartheta_{2}\right], \ldots,\left[\vartheta_{s}\right]\right\rangle$$ by \cite[Lemma 17]{hegazi2016classification}. 

\begin{thm}\cite[Theorem 18]{hegazi2016classification}\label{thm2}
	There exists a one-to-one correspondence between the set of $\mathrm{Aut}(A)$-orbits on $T_{s}(A)$ and the set of isomorphism classes of $E(A, V)$. This correspondence is defined by
	$$\operatorname{\mathrm{Orb}}\left\langle\left[\theta_{1}\right],\left[\theta_{2}\right], \ldots,\left[\theta_{s}\right]\right\rangle \in\left\{\operatorname{\mathrm{Orb}}(W): W \in T_{s}(A)\right\} \leftrightarrow\left[A_{\theta}\right] \in\left\{\left[A_{\vartheta}\right]: A_{\vartheta} \in E(A, V)\right\}$$ where $\theta(x, y)=\sum_{i=1}^{s} \theta_{i}(x, y) e_{i}$. 
\end{thm}

By this theorem, we may construct all anticommutative algebras of dimension $n$ with $s$-dimensional annihilator, given those algebras of dimension $n-s$, in the following way:
\begin{enumerate}
	\item For a given anticommutative algebra $A$ of dimension $n-s$, determine $H^{2}(A, \mathbb{C})$,  $\mathrm{\mathrm{Ann}}(A)$ and $\operatorname{\mathrm{Aut}}(A)$.
	\item Determine the set of $\mathrm{Aut}(A)$-orbits on $T_{s}(A)$.
	\item  For each orbit, construct the anticommutative algebra corresponding to a representative of it.
\end{enumerate}
This method gives all (Tortkara and non-Tortkara) anticommutative algebras with non-trivial annihilator. We want to develop this method in such a way that it only gives non-Tortkara anticommutative algebras. Clearly, any annihilator extension of non-Tortkara anticommutative algebra is non-Tortkara. So, we only have to study the central extensions of Tortkara algebras. Let $A$ be a Tortkara algebra and $\theta \in Z^{2}(A, \mathbb{C})$. Then $A_{\theta}$ is a Tortkara algebra if and only if $\theta((ab),(cb)) =\theta( J(a, b, c),b)$ for all $a, b, c \in A$, where $J(a, b, c)=(ab)c + (bc)a + (ca)b.$ Define a subspace $Z_{T}^{2}(A, \mathbb{C})$ of $Z^{2}(A, \mathbb{C})$ by
$$
Z_{T}^{2}(A, \mathbb{C})=\left\{\theta \in Z^{2}(A, \mathbb{C}): \theta((ab),(cb)) =\theta( J(a, b, c),b)\,\forall a,b,c \in A\right\}
$$
Define $H_{T}^{2}(A, \mathbb{C})=Z_{T}^{2}(A, \mathbb{C}) / B^{2}(A, \mathbb{C}) .$ Therefore, $H_{T}^{2}(A, \mathbb{C})$ is a subspace of $H^{2}(A, \mathbb{C})$. Define $R_{s}(A)= T_{s}(A)\cap G_{s}\left(H_{T}^{2}(A, \mathbb{C})\right)$. Then $T_{s}(A)=R_{s}(A) \cup U_{s}(A)$ where $U_{s}(A)=T_{s}(A)-R_{s}(A) .$ The sets $R_{s}(A)$ and $U_{s}(A)$ are stable under the action of $\mathrm{Aut}(A)$. Let $E_{T}(A, V)=\left\{A_{\theta} \in E(A, V): A_{\theta}  \mathrm{\,\,is\,\,Tortkara\,\, algebra} \right\} .$ Then $E(A, V)=E_{T}(A, V) \cup E_{non-T}(A, V)$ where $E_{ {non-T }}(A, V)=E(A, V)-E_{T}(A, V) .$

\begin{thm}\cite[Theorem 19]{hegazi2016classification}\label{thm3}
	Let $A$ be a Tortkara algebra.
	\begin{enumerate}
		\item There exists a one-to-one correspondence between the set of $\mathrm{Aut} (A)$-orbits on $R_{s}(A)$ and the set of isomorphism classes of $E_{T}(A, V)$.
		\item There exists a one-to-one correspondence between the set of $ \mathrm{Aut}(A)$-orbits on $U_{s}(A)$ and the set of isomorphism classes of $E_{n o n-T}(A, V)$.
	\end{enumerate}
\end{thm}

By this theorem and Theorem \ref{thm3}, we may construct all non-Tortkara anticommutative algebras of dimension $n$ with $s$-dimensional annihilator, given those algebras of dimension $n-s$, in the following way:
\begin{enumerate}
	\item 
	For a given anticommutative algebra $A$ of dimension $n-s$, if $A$ is non-Tortkara then do the following:
	\begin{enumerate}
		\item 
		Determine $H^{2}(A, \mathbb{C})$,  $\mathrm{\mathrm{\mathrm{Ann}}}(A)$ and  $\mathrm{Aut} (A)$.
		\item 
		Determine the set of $\mathrm{A u t}(A)$-orbits on $T_{s}(A)$.
		\item
		For each orbit, construct the anticommutative algebra corresponding to a representative of it.
	\end{enumerate}
	
	\item 
	Otherwise, do the following:
	\begin{enumerate}
		\item 
		Determine $H_{T}^{2}(A, \mathbb{C}), H^{2}(A, \mathbb{C}), \mathrm{A n n}(A)$ and $\mathrm{A u t}(A)$.
		\item 
		Determine the set of $\mathrm{Aut}(A)$-orbits on $U_{s}(A)$.
		\item
		For each orbit, construct the anticommutative algebra corresponding to a representative of it.
	\end{enumerate}
\end{enumerate}

\section{The proof of Theorem \ref{thm1}}\label{sec3}
Thanks to \cite{kaygorodov2020algebraic}, we have the classification of all nontrivial 4-dimensional nilpotent anticommutative algebras.
$$
\begin{array}{|l|l|l|l|}
	\hline {A} & \text {Multiplication table } & {H}_{{T}}^{2}({A},\mathbb{C}) & {H}^{2}({A},\mathbb{C}) \\
	
	\hline {A}_{1} & \rm{trivial\,\,algebra}  &
	\begin{aligned}
		&\left\langle\left[\Delta_{12}\right],\left[\Delta_{13}\right],\left[\Delta_{14}\right],\right.\\ &\left.\left[\Delta_{23}\right],\left[\Delta_{24}\right],\left[\Delta_{34}\right]\right\rangle
	\end{aligned} &
	\begin{aligned}
	&\left\langle\left[\Delta_{12}\right],\left[\Delta_{13}\right],\left[\Delta_{14}\right],\right.\\ &\left.\left[\Delta_{23}\right],\left[\Delta_{24}\right],\left[\Delta_{34}\right]\right\rangle
	\end{aligned}\\
	\hline {A}_{2} & e_{1} e_{2}=e_{3} & \left\langle\left[\Delta_{13}\right],\left[\Delta_{14}\right],\left[\Delta_{23}\right],\left[\Delta_{24}\right],\left[\Delta_{34}\right]\right\rangle & {H}_{{T}}^{2}\left({A}_{2},\mathbb{C}\right) \\
	\hline {A}_{3} & e_{1} e_{2}=e_{3}, e_{1} e_{3}=e_{4} & \left\langle\left[\Delta_{14}\right],\left[\Delta_{23}\right],\left[\Delta_{24}\right]\right\rangle & {H}_{{T}}^{2}\left({A}_{3},\mathbb{C}\right) \oplus\left\langle\left[\Delta_{34}\right]\right\rangle \\
	\hline
\end{array}
$$
In view of \cite{gorshkov2019variety}, all anticommutative central extensions of $A_{2}$ and of the 4-dimensional trivial algebra
are Tortkara algebras, so we need only consider central extensions of $A_{3}$. After direct calculation, we can get the form of $\mathrm{Aut}(A_3)$:
$$
\phi=\left(\begin{array}{cccc}
	x & 0 & 0 & 0 \\
	y & z & 0 & 0 \\
	u & v & x z & 0 \\
	h & g & x v & x^{2} z
\end{array}\right)
$$
where $xz\not=0$.

In the following subsections we will give the classification of all $n$-dimensional nilpotent non-Tortkara anticommutative algebras with $(n-4)$-dimensional annihilator.
\subsection{$n=4$}
There is no anticommutative algebra satisfies conditions of Theorem \ref{thm1}, because the annihilator of $A_1,A_2$ and $A_3$ are non-empty.

\subsection{$n=5$}\label{sec3.2}
The anticommutative algebra satisfying conditions of Theorem \ref{thm1} must have no annihilator component. It is non-split non-Tortkara central extensions of $A_3$. According to Theorem \ref{thm3}, we
need to find the representatives of the $\mathrm{Aut} (A_3)$-orbits on $U_1 (A_3)$. Choose an arbitrary subspace $W\in U_1(A_3)$. Such a subspace is spanned by $\theta=a_1[\Delta_{14}]+a_2[\Delta_{23}]+a_3[\Delta_{24}]+a_4[\Delta_{34}]$ with $a_4\not=0$. Let $\phi\in \mathrm{Aut}(A_3)$ be the following automorphism   
$$
\phi=\left(\begin{array}{cccc}
	x & 0 & 0 & 0 \\
	y & z & 0 & 0 \\
	-(a_1x+a_3y) & -a_3z & x z & 0 \\
	h & a_2z & -a_3zx  & x^{2} z
\end{array}\right)
$$
where $xz\not=0$. Then $\phi W=<[\Delta_{34}]>$. Hence we get a representative $<[\Delta_{34}]>$. This shows that $\mathrm{Aut} (A_3)$ has only one orbit on  $U_1 (A_3)$. So we get the algebra:
\[
A_{5,1}:e_{1} e_{2}=e_{3}, e_{1} e_{3}=e_{4}, e_{3} e_{4}=e_{5}.
\]
\subsection{$n=6$}\label{sec3.3}
First we classify nilpotent non-Tortkara anticommutative algebras with annihilator component. We get the algebras $A_{6,1}=A_{5,1}\oplus \mathbb{C}e_6$. Next we classify nilpotent non-Tortkara anticommutative without any annihilator component. For this, choose an arbitrary subspace $W\in U_2(A_3)$. Such a subspace is spanned by $[\theta_1]=a_1[\Delta_{14}]+a_2[\Delta_{23}]+a_3[\Delta_{24}]+a_4[\Delta_{34}]$ and $[\theta_2]=b_1[\Delta_{14}]+b_2[\Delta_{23}]+b_3[\Delta_{24}]+b_4[\Delta_{34}]$
with $(a_4,b_4)\not=(0,0)$. By possibly swapping $[\theta_1]$ and $[\theta_2]$, we may assume with out loss of generality that $a_4\not=0$. Then, from Subsection \ref{sec3.2}, we may assume that $[\theta_1]=[\Delta_{34}]$. Further, by subtracting scalar multiples of  $[\theta_1]$ from $[\theta_2]$, we may assume
that $[\theta_2]=a_1[\Delta_{14}]+a_2[\Delta_{23}]+a_3[\Delta_{24}]$. Hence we may assume without loss of generality that $$W=<[\Delta_{34}],a_1[\Delta_{14}]+a_2[\Delta_{23}]+a_3[\Delta_{24}]>.
$$
Let us consider the following cases:

\textbf{Case1:} $a_3\not=0$, Let $\phi\in \mathrm{Aut}(A_3)$ be as follows:
$$
\phi=\left(\begin{array}{cccc}
	x & 0 & 0 & 0 \\
	-a_1xa_3^{-1} & z & 0 & 0 \\
	0 & -a_2za_3^{-1}& x z & 0 \\
	h & v^2z^{-1} & -a_2za_3^{-1}x  & x^{2} z
\end{array}\right)
$$
where $xz\not=0$. Then $\phi W=<vzx^2[\Delta_{24}]+x^3z^2[\Delta_{34}],a_3x^2z^2[\Delta_{24}]>=<[\Delta_{34}],[\Delta_{24}]>=W_1$

\textbf{Case2:} $a_3=a_1=0,a_2\not=0$, Let $\phi\in \mathrm{Aut}(A_3)$ be as follows:
$$
\phi=\left(\begin{array}{cccc}
	x & 0 & 0 & 0 \\
	y & z & 0 & 0 \\
	0 & 0 & x z & 0 \\
	h & 0 & 0 & x^{2} z
\end{array}\right)
$$
where $xz\not=0$. Then $\phi W=<x^3z^2[\Delta_{34}],a_2xz^2[\Delta_{23}]>=<[\Delta_{34}],[\Delta_{23}]>=W_2$

\textbf{Case3:} $a_3=a_2=0,a_1\not=0$, Let $\phi\in \mathrm{Aut}(A_3)$ be as follows:
$$
\phi=\left(\begin{array}{cccc}
	x & 0 & 0 & 0 \\
	y & z & 0 & 0 \\
	0 & 0 & x z & 0 \\
	h & 0 & 0 & x^{2} z
\end{array}\right)
$$
where $xz\not=0$. Then $\phi W=<x^3z^2[\Delta_{34}],a_1zx^3[\Delta_{14}]>=<[\Delta_{34}],[\Delta_{14}]>=W_3$

\textbf{Case4:} $a_3=0,a_1a_2\not=0$, Let $\phi\in \mathrm{Aut}(A_3)$ be as follows:
$$
\phi=\left(\begin{array}{cccc}
	x & 0 & 0 & 0 \\
	0 & z & 0 & 0 \\
	0 & 0 & x z & 0 \\
	h & 0 & 0 & x^{2} z
\end{array}\right)
$$
where $x=a_2,z=a_1a_2$. Then $\phi W=<x^3z^2[\Delta_{34}],a_1^2a_2^4([\Delta_{23}]+[\Delta_{14}])>=<[\Delta_{34}],[\Delta_{23}]+[\Delta_{14}]>=W_4$

As shown we have three representatives, namely
\[
\begin{aligned}
	&W_1:<[\Delta_{34}],[\Delta_{24}]>\\
	&W_2:<[\Delta_{34}],[\Delta_{23}]>\\
	&W_3:<[\Delta_{34}],[\Delta_{14}]>\\
	&W_4:<[\Delta_{34}],[\Delta_{23}]+[\Delta_{14}]>\\
\end{aligned}
\]
Next, we claim that $\mathrm{Orb}(W_i)\cap \mathrm{Orb}(W_j)=\emptyset,i\not=j$. Choose any element of $\mathrm{Aut}(A_3)$:
$$
\phi=\left(\begin{array}{cccc}
	x & 0 & 0 & 0 \\
	y & z & 0 & 0 \\
	u & v & x z & 0 \\
	h & g & x v & x^{2} z
\end{array}\right)
$$
where $xz\not=0$. Then
\[
\phi W_1=<x^2z^2[\Delta_{34}]-zg[\Delta_{23}]+x(uz-yv)[\Delta_{14}],xz[\Delta_{24}]+v[\Delta_{23}]+xy[\Delta_{14}]>
\]
$[\Delta_{23}]$ is not in $\phi W_1$, otherwise there are $\lambda_1$ and $\lambda_2$ such that
\[
[\Delta_{23}]=\lambda_1(x^2z^2[\Delta_{34}]-zg[\Delta_{23}]+x(uz-yv)[\Delta_{14}])+\lambda_2(xz[\Delta_{24}]+v[\Delta_{23}]+xy[\Delta_{14}])
\]
We have $\lambda_1=\lambda_2=0$. Obviously, the above equation cannot be holded, it's a contradiction. Then $\mathrm{Orb}(W_1)\cap \mathrm{Orb}(W_2)=\emptyset$. $[\Delta_{14}]$ is not in $\phi W_1$, otherwise there are $\lambda_1$ and $\lambda_2$ such that
\[
[\Delta_{14}]=\lambda_1(x^2z^2[\Delta_{34}]-zg[\Delta_{23}]+x(uz-yv)[\Delta_{14}])+\lambda_2(xz[\Delta_{24}]+v[\Delta_{23}]+xy[\Delta_{14}])
\]
We have $\lambda_1=\lambda_2=0$. Obviously, the above equation cannot be holded, it's a contradiction. Then $\mathrm{Orb}(W_1)\cap \mathrm{Orb}(W_3)=\emptyset$. $[\Delta_{14}]+[\Delta_{23}]$ is not in $\phi W_1$, otherwise there are $\lambda_1$ and $\lambda_2$ such that
\[
[\Delta_{14}]+[\Delta_{23}]=\lambda_1(x^2z^2[\Delta_{34}]-zg[\Delta_{23}]+x(uz-yv)[\Delta_{14}])+\lambda_2(xz[\Delta_{24}]+v[\Delta_{23}]+xy[\Delta_{14}])
\]
We have $\lambda_1=\lambda_2=0$. Obviously, the above equation cannot be holded, it's a contradiction. Then $\mathrm{Orb}(W_1)\cap \mathrm{Orb}(W_4)=\emptyset$.

After $W_2$ is acted on by $\phi$, we can get
\[
\phi W_2=<x^2z^2[\Delta_{34}]+zxv[\Delta_{24}]+uzx[\Delta_{14}],[\Delta_{23}]>
\]
$[\Delta_{14}]$ is not in $\phi W_2$, otherwise there are $\lambda_1$ and $\lambda_2$ such that
\[
[\Delta_{14}]=\lambda_1(x^2z^2[\Delta_{34}]+zxv[\Delta_{24}]+uzx[\Delta_{14}])+\lambda_2[\Delta_{23}]
\]
We have $\lambda_1=\lambda_2=0$, Obviously, the above equation cannot be holded, it's a contradiction. Then $\mathrm{Orb}(W_2)\cap \mathrm{Orb}(W_3)=\emptyset$. $[\Delta_{14}]+[\Delta_{23}]$ is not in $\phi W_2$, otherwise there are $\lambda_1$ and $\lambda_2$ such that
\[
[\Delta_{14}]+[\Delta_{23}]=\lambda_1(x^2z^2[\Delta_{34}]+zxv[\Delta_{24}]+uzx[\Delta_{14}])+\lambda_2[\Delta_{23}]
\]
We have $\lambda_1=0,\lambda_2=1$. Obviously, the above equation cannot be holded, it's a contradiction. Then $\mathrm{Orb}(W_2)\cap \mathrm{Orb}(W_4)=\emptyset$. 

After $W_3$ is acted on by $\phi$, we can get
\[
\phi W_3=<x^2z^2[\Delta_{34}]+zxv[\Delta_{24}]+(v^2-zg)[\Delta_{23}],[\Delta_{14}]>
\]
$[\Delta_{14}]+[\Delta_{23}]$ is not in $\phi W_3$, otherwise there are $\lambda_1$ and $\lambda_2$ such that
\[
[\Delta_{14}]+[\Delta_{23}]=\lambda_1(x^2z^2[\Delta_{34}]+zxv[\Delta_{24}]+(v^2-zg)[\Delta_{23}])+\lambda_2[\Delta_{14}]
\]
We have $\lambda_1=0,\lambda_2=1$. Obviously,  the above equation cannot be holded, it's a contradiction. Then $\mathrm{Orb}(W_3)\cap \mathrm{Orb}(W_4)=\emptyset$. All in all, we have $\mathrm{Orb}(W_i)\cap \mathrm{Orb}(W_j)=\emptyset,i\not=j$. So, the corresponding anticommutative algebras are pairwise non-isomorphic. Hence
we get the algebras:
\[
\begin{aligned}
	&A_{6,2}:e_{1} e_{2}=e_{3}, e_{1} e_{3}=e_{4}, e_{3} e_{4}=e_{5}, e_{2} e_{4}=e_{6};\\
	&A_{6,3}:e_{1} e_{2}=e_{3}, e_{1} e_{3}=e_{4}, e_{3} e_{4}=e_{5}, e_{2} e_{3}=e_{6};\\
	&A_{6,4}:e_{1} e_{2}=e_{3}, e_{1} e_{3}=e_{4}, e_{3} e_{4}=e_{5}, e_{1} e_{4}=e_{6};\\
	&A_{6,5}:e_{1} e_{2}=e_{3}, e_{1} e_{3}=e_{4}, e_{3} e_{4}=e_{5}, e_{2} e_{3}=e_{6}, e_{1} e_{4}=e_{6}.\\
\end{aligned}
\]

\subsection{$n=7$}
First we classify the non-Tortkara anticommutative algebras with annihilator component. We get the algebras $A_{7,i}=A_{6,i}\oplus\mathbb{C}e_7,i=1,\cdots,5$. Next we classify the non-Tortkara anticommutative algebras without any annihilator component. Choose an arbitrary subspace
$W\in U_3(A_3)$. Such a subspace is spanned by
\[
\begin{aligned}
	&[\theta_1]=a_1[\Delta_{14}]+a_2[\Delta_{23}]+a_3[\Delta_{24}]+a_4[\Delta_{34}], \\
	 &[\theta_2]=b_1[\Delta_{14}]+b_2[\Delta_{23}]+b_3[\Delta_{24}]+b_4[\Delta_{34}]\\
	&[\theta_3]=c_1[\Delta_{14}]+c_2[\Delta_{23}]+c_3[\Delta_{24}]+c_4[\Delta_{34}]
\end{aligned}
\]
with $(a_4,b_4,c_4)\not=(0,0,0)$.  In view of Subsection \ref{sec3.3}, we may assume without loss of generality that $W\in\{S_1, S_2, S_3,S_4\}$ where
$$
\begin{aligned}
	&S_1=<[\Delta_{34}],[\Delta_{24}],a_1[\Delta_{14}]+a_2[\Delta_{23}]>\\
	&S_2=<[\Delta_{34}],[\Delta_{14}],a_2[\Delta_{23}]+a_3[\Delta_{24}]>\\
	&S_3=<[\Delta_{34}],[\Delta_{23}],a_1[\Delta_{14}]+a_3[\Delta_{24}]>\\
	&S_4=<[\Delta_{34}],[\Delta_{14}]+[\Delta_{23}],a_1[\Delta_{14}]+a_3[\Delta_{24}]>
\end{aligned}
$$

\textbf{Case1:} $W=S_1$

\textbf{Case1.1:} $a_1=0$ or $a_2=0$, then $W=<[\Delta_{34}],[\Delta_{24}],[\Delta_{23}]>$ or $W=<[\Delta_{34}],[\Delta_{24}],[\Delta_{14}]>$.

\textbf{Case1.2:} $a_1a_2\not=0$. Let $\phi\in \mathrm{Aut}(A_3)$ be as follows:
$$
\phi=\left(\begin{array}{cccc}
	x & 0 & 0 & 0 \\
	0 & z & 0 & 0 \\
	0 & 0 & x z & 0 \\
	h & 0 & 0 & x^{2} z
\end{array}\right)
$$
where $x=a_2,z=a_1a_2$. Then $\phi W=<x^3z^2[\Delta_{34}],x^2z^2[\Delta_{24}],a_1^2a_2^4([\Delta_{14}]+[\Delta_{23}])>=<[\Delta_{34}],[\Delta_{24}],[\Delta_{14}]+[\Delta_{23}]>$

\textbf{Case2:} $W=S_2$

\textbf{Case2.1:} $a_2=0$ or $a_3=0$. Then $W=<[\Delta_{34}],[\Delta_{14}],[\Delta_{24}]>$ or $W=<[\Delta_{34}],[\Delta_{14}],[\Delta_{23}]$.

\textbf{Case2.2:} $a_2a_3\not=0$. Let $\phi\in \mathrm{Aut}(A_3)$ be as follows:
$$
\phi=\left(\begin{array}{cccc}
	x & 0 & 0 & 0 \\
	0 & z & 0 & 0 \\
	0 & 0 & x z & 0 \\
	h & 0 & 0 & x^{2} z
\end{array}\right)
$$
where $x=a_3^2,z=a_2a_3^{-1}$. Then $\phi W=<x^3z^2[\Delta_{34}],zx^3[\Delta_{14}],a_2^3([\Delta_{23}]+[\Delta_{24}])> =<[\Delta_{34}],[\Delta_{14}],[\Delta_{23}]+[\Delta_{24}]>$

\textbf{Case3:} $W=S_3$

\textbf{Case3.1:} $a_1=0$ or $a_3=0$. Then $W=<[\Delta_{34}],[\Delta_{23}],[\Delta_{24}]>$ or $W=<[\Delta_{34}],[\Delta_{23}],[\Delta_{14}]>$.

\textbf{Case3.2:} $a_1a_3\not=0$. Let $\phi\in \mathrm{Aut}(A_3)$ be as follows:
$$
\phi=\left(\begin{array}{cccc}
	x & 0 & 0 & 0 \\
	0 & z & 0 & 0 \\
	0 & 0 & x z & 0 \\
	h & 0 & 0 & x^{2} z
\end{array}\right)
$$
where $x=a_3,z=a_1$. Then $\phi W=<x^3z^2[\Delta_{34}],xz^2[\Delta_{23}],a_1^2a_3^3([\Delta_{14}]+[\Delta_{24}])> =<[\Delta_{34}],[\Delta_{23}],[\Delta_{14}]+[\Delta_{24}]>$.

\textbf{Case4:} $W=S_4$

\textbf{Case4.1:} $a_1=0$ or $a_3=0$. Then $W=<[\Delta_{34}],[\Delta_{14}]+[\Delta_{23}],[\Delta_{24}]>$ or $W=<[\Delta_{34}],[\Delta_{23}],[\Delta_{14}]>$.

\textbf{Case4.2:} $a_1a_3\not=0$. Let $\phi\in \mathrm{Aut}(A_3)$ be as follows:
$$
\phi=\left(\begin{array}{cccc}
	1 & 0 & 0 & 0 \\
	-a_1a_3^{-1}& 1 & 0 & 0 \\
	0 & 0 & 1 & 0 \\
	0 & 0 & 0 & 1
\end{array}\right)
$$
Then $\phi W=<[\Delta_{34}],[\Delta_{14}]+[\Delta_{23}],a_3[\Delta_{24}]>=<[\Delta_{34}],[\Delta_{14}]+[\Delta_{23}],[\Delta_{24}]>$

To summarize, we have the following representatives:
\[
\begin{aligned}
	&W_1=<[\Delta_{34}],[\Delta_{24}],[\Delta_{23}]>\\
	&W_2=<[\Delta_{34}],[\Delta_{24}],[\Delta_{14}]>\\
	&W_3=<[\Delta_{34}],[\Delta_{14}],[\Delta_{23}]>\\
	&W_4=<[\Delta_{34}],[\Delta_{24}],[\Delta_{14}]+[\Delta_{23}]>\\
	&W_5=<[\Delta_{34}],[\Delta_{14}],[\Delta_{23}]+[\Delta_{24}]>\\
	&W_6=<[\Delta_{34}],[\Delta_{23}],[\Delta_{14}]+[\Delta_{24}]>
\end{aligned}
\]

Let us now determine the possible orbits among the representatives $W_1,\cdots,W_6$. Choose any element of $\mathrm{Aut}(A_3)$:
$$
\phi=\left(\begin{array}{cccc}
	x & 0 & 0 & 0 \\
	y & z & 0 & 0 \\
	u & v & x z & 0 \\
	h & g & x v & x^{2} z
\end{array}\right)
$$
where $xz\not=0$. Then 
\[
\phi W_1=<[\Delta_{23}],z[\Delta_{24}]+y[\Delta_{14}],xz[\Delta_{34}]+\left(u-vyz^{-1}\right)[\Delta_{14}]>
\]
Let $y=z,u=v$, we have $\phi W_1=W_6$. Then $\mathrm{Orb}(W_1)=\mathrm{Orb}(W_6)$. $[\Delta_{14}]$ is not in $\phi W_1$, otherwise there are $\lambda_1$ , $\lambda_2$ and $\lambda_3$  such that
\[
[\Delta_{14}]=\lambda_1[\Delta_{23}]+\lambda_2(z[\Delta_{24}]+y[\Delta_{14}])+\lambda_3(xz[\Delta_{34}]+\left(u-vyz^{-1}\right)[\Delta_{14}])
\]
We have $\lambda_1=\lambda_2=\lambda_3=0$. Obviously,  the above equation cannot be holded, it's a contradiction. Then $\mathrm{Orb}(W_1)\cap \mathrm{Orb}(W_2)=\emptyset$, $\mathrm{Orb}(W_1)\cap \mathrm{Orb}(W_3)=\emptyset$ and $\mathrm{Orb}(W_1)\cap \mathrm{Orb}(W_5)=\emptyset$. We claim that $\mathrm{Orb}(W_1)\cap \mathrm{Orb}(W_4)=\emptyset$, otherwise there are $\lambda_1$ , $\lambda_2$ and $\lambda_3$  such that
\[
[\Delta_{34}]=\lambda_1[\Delta_{23}]+\lambda_2(z[\Delta_{24}]+y[\Delta_{14}])+\lambda_3(xz[\Delta_{34}]+\left(u-vyz^{-1}\right)[\Delta_{14}])
\]
We have $\lambda_1=\lambda_2=0,u=vyz^{-1}$. There are $\lambda_1$ , $\lambda_2$ and $\lambda_3$  such that
\[
[\Delta_{24}]=\lambda_1[\Delta_{23}]+\lambda_2(z[\Delta_{24}]+y[\Delta_{14}])+\lambda_3(xz[\Delta_{34}]+\left(u-vyz^{-1}\right)[\Delta_{14}])
\]
We have $\lambda_1=\lambda_3=0,y=0$ and $u=y=0$. There are $\lambda_1$ , $\lambda_2$ and $\lambda_3$  such that
\[
[\Delta_{14}]+[\Delta_{23}]=\lambda_1[\Delta_{23}]+\lambda_2z[\Delta_{24}]+\lambda_3xz[\Delta_{34}]
\]
We have $\lambda_1=1,\lambda_2=\lambda_3=0$. Obviously,  the above equation cannot be holded, it's a contradiction. Then $\mathrm{Orb}(W_1)\cap \mathrm{Orb}(W_4)=\emptyset$.

After $W_2$ is acted on by $\phi$, we can get
\[
\phi W_2=<[\Delta_{14}],xz[\Delta_{24}]+v[\Delta_{23}],x^2z[\Delta_{34}]-g[\Delta_{23}]>
\]
Let $g=0,v=xz$, we have $\phi W_2=W_5$. Then $\mathrm{Orb}(W_2)=\mathrm{Orb}(W_5)$. We claim that $\mathrm{Orb}(W_2)\cap \mathrm{Orb}(W_3)=\emptyset$, otherwise there are $\lambda_1$ , $\lambda_2$ and $\lambda_3$  such that
\[
[\Delta_{23}]=\lambda_1[\Delta_{14}]+\lambda_2(xz[\Delta_{24}]+v[\Delta_{23}])+\lambda_3(x^2z[\Delta_{34}]-g[\Delta_{23}])
\]
We have $\lambda_1=\lambda_2=\lambda_3=0$. Obviously,  the above equation cannot be holded, it's a contradiction. Then $\mathrm{Orb}(W_2)\cap \mathrm{Orb}(W_3)=\emptyset$.

We claim that $\mathrm{Orb}(W_2)\cap \mathrm{Orb}(W_4)=\emptyset$, otherwise there are $\lambda_1$ , $\lambda_2$ and $\lambda_3$  such that
\[
[\Delta_{34}]=\lambda_1[\Delta_{14}]+\lambda_2(xz[\Delta_{24}]+v[\Delta_{23}])+\lambda_3(x^2z[\Delta_{34}]-g[\Delta_{23}])
\]
We have $g=0$. There are $\lambda_1$ , $\lambda_2$ and $\lambda_3$  such that
\[
[\Delta_{24}]=\lambda_1[\Delta_{14}]+\lambda_2(xz[\Delta_{24}]+v[\Delta_{23}])+\lambda_3(x^2z[\Delta_{34}]-g[\Delta_{23}])
\]
We have $v=0$. There are $\lambda_1$ , $\lambda_2$ and $\lambda_3$  such that
\[
[\Delta_{14}]+[\Delta_{23}]=\lambda_1[\Delta_{14}]+\lambda_2 xz[\Delta_{24}]+\lambda_3x^2z[\Delta_{34}]
\]
We have $\lambda_1=\lambda_2=\lambda_3=0$. Obviously,  the above equation cannot be holded, it's a contradiction. Then $\mathrm{Orb}(W_2)\cap \mathrm{Orb}(W_4)=\emptyset$.

After $W_3$ is acted on by $\phi$, we can get
\[
\phi W_3=<[\Delta_{14}],[\Delta_{23}],xz[\Delta_{34}]+v[\Delta_{24}]>
\]
We claim that $\mathrm{Orb}(W_3)\cap \mathrm{Orb}(W_4)=\emptyset$, otherwise there are $\lambda_1$ , $\lambda_2$ and $\lambda_3$  such that
\[
[\Delta_{34}]=\lambda_1[\Delta_{14}]+\lambda_2[\Delta_{23}]+\lambda_3(xz[\Delta_{34}]+v[\Delta_{24}])
\]
We have $v=0$. There are $\lambda_1$ , $\lambda_2$ and $\lambda_3$  such that
\[
[\Delta_{24}]=\lambda_1[\Delta_{14}]+\lambda_2[\Delta_{23}]+\lambda_3xz[\Delta_{34}]
\]
We have $\lambda_1=\lambda_2=\lambda_3=0$. Obviously,  the above equation cannot be holded, it's a contradiction. Then $\mathrm{Orb}(W_3)\cap \mathrm{Orb}(W_4)=\emptyset$.

All representatives of $\mathrm{Aut}(A_3)$-orbits are:
\[
\begin{aligned}
	&W_1=<[\Delta_{34}],[\Delta_{24}],[\Delta_{23}]>\\
	&W_2=<[\Delta_{34}],[\Delta_{24}],[\Delta_{14}]>\\
	&W_3=<[\Delta_{34}],[\Delta_{14}],[\Delta_{23}]>\\
	&W_4=<[\Delta_{34}],[\Delta_{24}],[\Delta_{14}]+[\Delta_{23}]>
\end{aligned}
\]
So the corresponding anticommutative algebras are pairwise non-isomorphic. Thus we get the algebras:
\[
\begin{aligned}
	&A_{7,6}:e_{1} e_{2}=e_{3}, e_{1} e_{3}=e_{4}, e_{3} e_{4}=e_{5}, e_{2} e_{4}=e_{6}, e_{2} e_{3}=e_{7};\\
	&A_{7,7}:e_{1} e_{2}=e_{3}, e_{1} e_{3}=e_{4}, e_{3} e_{4}=e_{5}, e_{2} e_{4}=e_{6}, e_{1} e_{4}=e_{7};\\
	&A_{7,8}:e_{1} e_{2}=e_{3}, e_{1} e_{3}=e_{4}, e_{3} e_{4}=e_{5}, e_{1} e_{4}=e_{6}, e_{2} e_{3}=e_{7};\\
	&A_{7,9}:e_{1} e_{2}=e_{3}, e_{1} e_{3}=e_{4}, e_{3} e_{4}=e_{5}, e_{2} e_{4}=e_{6}, e_{1} e_{4}=e_{7}, e_{2} e_{3}=e_{7}.
\end{aligned}
\]

\subsection{$n=8$}\label{sec3.5}
First we classify the non-Tortkara anticommutative algebras with annihilator component. We get the algebras $A_{8,i}=A_{7,i}\oplus\mathbb{C}e_8,i=1,\cdots,9$. Next we classify the non-Tortkara anticommutative algebras without any annihilator component. Choose an arbitrary subspace $W\in  U_4(A_3)$. Then $W=H^2(A_3,\mathbb{C})$. So we have only one orbit with a
representative $<[\Delta_{14}],[\Delta_{23}],[\Delta_{24}],[\Delta_{34}]>$. Hence we get the algebra:
\[
A_{8,10}:e_{1} e_{2}=e_{3}, e_{1} e_{3}=e_{4}, e_{3} e_{4}=e_{5}, e_{2} e_{4}=e_{6}, e_{2} e_{3}=e_{7},e_{1} e_{4}=e_{8}.
\]

\subsection{$n\geqslant 9$}
From the results of Subsection \ref{sec3.5}, we get the algebras:
$$
A_{n, i}=A_{8, i} \oplus \mathbb{C} e_{9} \oplus \cdots \oplus \mathbb{C} e_{n}, i=1, \ldots, 10.
$$
\section{Acknowledgment}
 I would like to express my sincere gratitude to professor Jianzhi Han from Tongji University.
\bibliographystyle{plain}  
\bibliography{refs}

\end{document}